\documentclass[11pt]{article}

\usepackage{amsmath, amssymb, amsthm}
\usepackage{geometry}
\usepackage{xcolor}
\usepackage[font=small,labelfont=bf, width=.85\textwidth]{caption}
\usepackage{enumerate}
\usepackage{authblk}
\usepackage{comment}
\usepackage{soul}

\newtheorem{theorem}{Theorem}

\newtheorem{conjecture}[theorem]{Conjecture}
\newtheorem*{conjecture*}{Conjecture}
\newtheorem{observation}[theorem]{Observation}

\newtheorem{definition}[theorem]{Definition}

\title{Balanced chromatic number and Hadwiger-like conjectures}

\author[1,5]{Andrea Jim\'enez} 
\author[2]{Jessica McDonald} 
\author[3]{Reza Naserasr}
\author[4]{Kathryn Nurse}
\author[1]{Daniel A. Quiroz}

\affil[1]{\small Instituto de Ingenier\'ia Matem\'atica-CIMFAV, Universidad de Valpara\'iso, Valpara\'iso, Chile. {Emails: \texttt{\{andrea.jimenez, daniel.quiroz\}@uv.cl}}}
\affil[2] {\small Auburn University, Auburn, AL, USA 36849. {Email: \texttt{mcdonald@auburn.edu}}}
\affil[3]{\small Université Paris Cité, CNRS, IRIF, F-75013, Paris, France. {Email: \texttt{reza@irif.fr}}}	
\affil[4]{\small \'Ecole Normale Sup\'erieure, Paris, France. {Email: \texttt{kathryn.bale@gmail.com}}}
\affil[5]{Millennium Nucleus for Social Data Science (SODAS), Santiago, Chile}

\begin{document}

	\date{}

	\maketitle
	
	\begin{abstract} Motivated by different characterizations of planar graphs and the 4-Color Theorem, several structural results concerning graphs of high chromatic number have been obtained. Toward strengthening some of these results, we consider the \emph{balanced chromatic number}, $\chi_b(\hat{G})$, of a signed graph $\hat{G}$. This is the minimum number of parts into which the vertices of a signed graph can be partitioned so that none of the parts induces a negative cycle. This extends the notion of the chromatic number of a graph since $\chi(G)=\chi_b(\tilde{G})$, where $\tilde{G}$ denotes the signed graph obtained from~$G$ by replacing each edge with a pair of (parallel) positive and negative edges. We introduce a signed version of Hadwiger's conjecture as follows.
		
		\begin{conjecture*}			
			If a signed graph $\hat{G}$ has no negative loop and no $\tilde{K}_t$-minor, then its balanced chromatic number is at most $t-1$.
		\end{conjecture*}
		
		We prove that this conjecture is, in fact, equivalent to  Hadwiger's conjecture and show its relation to the odd Hadwiger Conjecture. 

		Motivated by these results, we also consider the relation between subdivisions and balanced chromatic number. We prove that if $(G, \sigma)$ has no negative loop and no $\tilde{K}_t$-subdivision, then it admits a balanced $\frac{79}{2}t^2$-coloring. This qualitatively generalizes a result of Kawarabayashi (2013) on totally odd subdivisions.
		
		Finally, following supportive results in the literature on the fractional variant of Hadwiger's conjecture, we show that the fractional balanced chromatic number of any signed graph with no positive loop and no $\tilde{K}_t$-minor is at most $2t-2$.
		
	\end{abstract}
	
	\section{Introduction}
	
	The 4-Color Theorem (4CT) has been the driving engine behind many of the developments in graph theory. The characterization of planar graphs as the class of graphs with no $K_5$ or $K_{3,3}$-minor (Wagner~\cite{W1937}), or as the class of graphs with no $K_5$- or $K_{3,3}$-subdivision (Kuratowski \cite{K1930}) has led to various conjectures which generalize the 4CT, mostly in the form of the following question: Given a graph of high chromatic number, what sort of structure(s) are we sure to have in $G$? A stronger version of the 4CT, obtained through Wagner's characterization of $K_5$-minor-free graphs \cite{W1937},
	tells us that if the chromatic number of $G$ is at least 5, then $G$ must contain a $K_5$-minor. 
	On the other hand, the 4CT and Kuratowski's characterization of planar graphs 
	imply that every 5-chromatic graph contains either a $K_5$- or a $K_{3,3}$-subdivision, and it remains open whether we can always find a $K_5$-subdivision in such a graph. 
	To consider similar problems for graphs with arbitrary chromatic number we use the following notation.
	
	\begin{definition}
		Given a positive integer $t$, let $\mathcal{C}_t$ be the class of graphs of chromatic number at least $t$. Then $M(t)$ and $T(t)$ are defined to be, respectively, the largest $k$ and $\ell$ such that every member of $\mathcal{C}_t$ contains $K_k$ as a minor and $K_\ell$ as a subdivision.
	\end{definition}
	
	Using this terminology, the well-known Hadwiger's conjecture is reformulated as follows.

	\begin{conjecture} [Hadwiger~\cite{H1943}]
		$M(t)=t$ for every positive integer $t$.
	\end{conjecture}
	
	Similarly, H\'ajos conjectured that $T(t)=t$. However, H\'ajos' conjecture has been shown to be false for $t\geq 7$, remains open for $t=5,6$, and is verified for $t=3,4$. We refer the reader to \cite{XXYY21} and the references therein for more on the open cases of H\'ajos' conjecture. 
	
	Currently the best general upper bound known on the chromatic number of $K_t$-minor-free graphs is that of Delcourt and Postle \cite{DP21} who prove a $O(t\log \log t)$ upper bound. The best known general upper bound for the chromatic number of graphs with no  $K_t$-subdivision is $O(t^2)$ proved independently in \cite{BT98} and \cite{KS96}.

	While these results 
	are among the most famous results  on the structure of graphs of high chromatic number, their strength (and even the strength of the conjectured values, if proved to be correct) has also been subject to challenge. To best express this let us first look at the case $t=3$ of Hadwiger's conjecture: the set of $K_3$-minor-free graphs is that of acyclic graphs. While it is true that they are 2-colorable, we readily have a much stronger result: a graph is 2-colorable if (and only if) it has no odd cycle. More generally, given a graph $G$ if even just one edge is blown up to a large complete bipartite graph, then while the chromatic number of the resulting graph $G'$ remains the same, the order of clique-minor or  clique-subdivision one can find in $G'$ increases with the order of the complete bipartite graph considered.

	This issue has been dealt with in the literature by refining the containment relations in question in several different ways. The starting point is a result of Catlin \cite{C1979} who proved that if $G$ is not 3-colorable then it has a subdivision of $K_4$, the  triangles of which correspond to odd cycles in $G$. From here, two extensions have been proposed. On the one hand,  Gerards and Seymour, independently, introduced the odd Hadwiger Conjecture (presented in Section~\ref{sec:sign}) which, as we will see, is rooted in the idea of minors in signed graphs. On the other hand, further inspired by a result of Zang \cite{Z1998}, various extensions of results for subdivisions in graphs of high chromatic number have been obtained. 
	
	While these extensions manage to deal with the above-mentioned issues, proving results towards the corresponding conjectures can be particularly complicated. For instance, when working towards the odd Hadwiger Conjecture, induction cannot be easily applied because 
    based on the definition of minor for signed graphs given in the next section in order to apply the contraction operation one must first apply a switching which takes one outside of the desired class of (signed) graphs. Working within the more general framework of signed graphs can help in this regard and in this paper we aim to identify notions of coloring in signed graphs which are most suitable to deal with such problems. 
	
	The rest of the paper is organized as follows. In Section~\ref{sec:sign} we discuss signed graphs, their minors and their colorings. In Section \ref{sec:ourHad} we present a conjecture on balanced colorings of signed graphs that trivially generalizes Hadwiger's conjecture but is, in fact, equivalent to it. We also study our conjecture's relation to the odd Hadwiger's conjecture. Finally, in Section \ref{sec:sudsigned} we study the presence of subdivisions in signed graphs with high balanced chromatic number, qualitatively generalizing a result of Kawarabayashi \cite{K2013}.

	\section{Signed graphs}\label{sec:sign}
	Signed graphs offer a more complete model of networks (such as social ones), as compared to graphs. While a graph model for a network can only capture if two objects of the network are joined or not, in a signed graph model such a connection can be of two possible types: positive and negative.

	We use $(G, \sigma)$ to denote a \emph{signed graph} where $\sigma(e)$ determines the sign of the edge $e$ in $G$. When the \emph{signature} $\sigma$ makes no difference, we take $\hat{G}$ as a signed graph. The signed graph on $G$ where all edges are negative is denoted by $(G, -)$. For the purpose of this work, the most natural interpretation of graphs as a subclass of signed graphs is to see them as the class of signed graphs with all edges being negative.

	Given a graph $G$, the signed graph $\tilde{G}$ is built from $G$ by replacing each edge with a \emph{digon} which consists of two (parallel) edges: one positive and one negative. Signed graphs are normally equipped with the following basic but key mathematical operation: given a vertex $v$, a \emph{switching at $v$} consists of multiplying the signs of all edges incident to $v$ by a $-$. The \emph{sign} of a substructure (subgraph, minor, etc.) of $(G, \sigma)$ is the product of the signs of the edges of such structure (considering multiplicity). A key observation is that signs of cycles and closed walks are invariant under the switching operation.
	
	Signed graphs in this paper are permitted to have loops and parallel edges, unless otherwise stated. When stating results about coloring however, negative loops will never be considered.  
	
	\subsection{Minors of signed graphs}
	
	The notion of minors for signed graphs mirrors natural grouping in a social network: One can be added to a group if one already has a ``positive" relation with some member of the group. Formally, a signed graph $(H, \pi)$ is said to be a \emph{minor} of the signed graph $(G, \sigma)$ if it is obtained from $(G, \sigma)$ by a series of vertex  or edge deletions, contraction of positive edges, and switchings.  Just as with switching, the contraction operation does not change the sign of a cycle. Thus, unless a cycle is deleted, its image in $(H, \pi)$ is a closed walk of the same sign. Within graphs (noting that sign takes the role of parity in unsigned graphs) this means that the parity of signs is preserved, and this allows for graphs excluding $K_3$ as an odd-minor to be precisely the graphs with no odd cycle.
	
	The odd Hadwiger Conjecture, proposed independently by Gerards and Seymour is the following strengthening of Hadwiger's conjecture. We refer to \cite{JT95} (page 115) for the first citation and to \cite{NRS15} for this formulation in the language of signed graphs.
	
	\begin{conjecture}[odd Hadwiger]\label{conj:oddHadwiger}
		If $\chi(G)\geq t$, then $(G,-)$ has a $(K_t,-)$-minor.
	\end{conjecture}
	
	As a general signature is not needed in this statement, the conjecture is normally presented in the language of 2-vertex colored graphs where the colors actually indicate whether a switching has occurred at a vertex or not. For references and some earlier work on the conjecture see \cite{GGRSV06}.
	
	We define $\mathcal{OC}_{t}$ to be the class of signed graphs $(G,-)$ where $\chi(G)\geq t$. Let $OM(t)$ be the largest $k$ such that $(K_k,-)$ is a minor of every element of $\mathcal{OC}_{t}$. In this language, the odd Hadwiger Conjecture claims that $OM(t)=t$. 
	
	\subsection{Coloring signed graphs}
	A signed graph is said to be \emph{balanced} if it contains no negative cycles. This is equivalent to finding an edge-cut whose edges are all negative with all other edges positive (see \cite{H53} and \cite{Z82b} for a generalization). A subset $X$ of vertices of a signed graph $\hat{G}$ is said to be balanced if it induces no negative cycle. In other words, if the induced subgraph is balanced. 
	
	\begin{definition}
		A \emph{balanced $k$-coloring} of a signed graph $(G, \sigma)$ is to partition (or cover) its vertices with $k$ (possibly empty) balanced sets. The \emph{balanced chromatic number} of a signed graph $\hat{G}$, denoted $\chi_b(\hat{G})$, is the minimum $k$ for which $\hat{G}$ admits a balanced $k$-coloring.
	\end{definition}	
	
	It is obvious from the definition that $\hat{G}$ admits a balanced coloring for some $k$ ($k\leq |V(\hat{G})|$) if and only if it has no negative loops. If it has a negative loop then we may write $\chi_b(\hat{G})=\infty$. The notion of balanced coloring generalizes the notion of proper coloring of graphs by the observation that $\chi(G)=\chi_b(\tilde{G})$.
	
	The notion first appears in \cite{Z1987} where the term \emph{balanced partition number} is used. In practice we work with covering vertex set rather than partitioning it. Furthermore to be inline with main stream graph theory, we rather use the terms balanced coloring and balanced chromatic number. A closely related notion to balanced coloring is the notion of 0-free coloring also introduced by Zaslavsky in \cite{Z1982a}: a 0-free $k$-coloring of a signed graph $(G, \sigma)$ with no positive loop is an assignment $\psi$ of values from $\{\pm 1, \pm 2, \ldots, \pm k\}$ to the vertices of $G$ such that for each edge $e=uv$ we have $\psi(u)\neq \sigma(e)\psi(v)$. It is not difficult to observe that $(G, \sigma)$ admits a balanced $k$-coloring if and only if $(G, -\sigma)$ admits a 0-free $k$-coloring. 
	Thus some results on $0$-free colorings apply to balanced colorings as well. For example, it follows from a result of \cite{MRS16} that every $2k$-degenerate signed simple graph admits a balanced $k$-coloring.
	
	A refinement of balanced coloring is the notion of circular coloring of signed graphs, first given in \cite{NWZ21}. As in the case of $0$-free coloring, the definition is slightly modified here to better suit the relation with minor theory. A \emph{circular $r$-coloring} ($r\geq 2$) of a signed graph $\hat{G}$ is an assignment $\phi$ of the vertices of $\hat{G}$ to the points of a circle $O$ of circumference $r$ such that for each negative edge $xy$, $\phi(x)$ and $\phi(y)$ are at distance at least $1$, and that for each positive edge $zt$, $\phi(z)$ and $\phi(t)$ are at distance at most $\frac{r}{2}-1$ (equivalently, the distance of $\phi(z)$ from the antipodal of $\phi(t)$ is at least 1). Given a signed graph $\hat{G}$ with no negative loops, the smallest $r$ for which $\hat{G}$ admits an $r$-coloring is called the \emph{circular chromatic number} of $\hat{G}$, denoted $\chi_{c}(\hat{G})$. 
	
	The circular chromatic number of signed graphs is a refinement of the balanced chromatic number of signed graphs by the following formula whose proof is based on the notions of homomorphism given in the next section.
	
	$$\chi_b(\hat{G})=\Bigl\lceil \frac{\chi_{c}(\hat{G})}{2} \Bigr\rceil .$$ 
	
	\subsection{Homomorphisms of signed graphs}
	
	A \emph{homomorphism} of a signed graph $(G,\sigma)$ to a signed graph $(H,\pi)$ is a mapping of vertices and edges of $G$ to  the vertices and edges of $H$ (respectively) such that incidences, adjacencies, and signs of closed walks are all preserved. The condition of preserving signs of closed walks is equivalent to having a homomorphism from some switching equivalent version $(G, \sigma')$ to $(H,\pi)$ such that signs of the edges are preserved. For more on homomorphisms we refer to \cite{NSZ21}.
	
	Let $\tilde{K_k}^+$ be the signed graph on $k$ vertices where there is a positive loop on each vertex and any pair of distinct vertices are adjacent by both positive and negative edges. It is easily observed that:
	
	\begin{observation}
		The balanced chromatic number of a signed graph $(G, \sigma)$ (with no negative loop) is the smallest $k$ such that $(G, \sigma)$ admits a homomorphism to $\tilde{K_k}^+$.  
	\end{observation}
	
	Homomorphism targets for circular $r$-coloring of $(G,\sigma)$ are provided in \cite{NWZ21}. It turns out that $\tilde{K_k}^+$ is the homomorphism target for a signed graph to admit a circular $2k$-coloring. This proves the identity given in the previous section.
	
	\subsection{Connecting minor and homomorphisms of signed graphs} 
	
	The following theorem is the essential element for reducing problems on balanced chromatic number of minor closed families of signed graphs to problems of similar nature in graph (not signed).
	
	\begin{theorem}\label{thm:MinorAndHomImage}
		Given a signed graph $\hat{G}$ with no negative loop, there exists a signed graph $\tilde{H}$ with no negative loop but possibly with positive loops, where any pair of distinct connected vertices are connected by a digon, and such that $\tilde{H}$ is both a minor of $\hat{G}$ and a homomorphic image of it.
	\end{theorem}
	
	\begin{proof}
		Let $V_1$ be a maximal set of vertices of $(G,\sigma)$ that induces a connected balanced subgraph. Since it induces a balanced subgraph, there is a suitable switching on subset of vertices in $V_1$ after which all edges induced by $V_1$ are positive. After such a switching, identifying all vertices of $V_1$ into a single vertex, say $v_1$, and possibly adding a positive loop on it if needed, is both a homomorphism operation and a minor operation (the latter because $V_1$ induces a connected subgraph).
		Let $\hat{G_1}$ be the resulting signed graph. In $\hat{G_{1}}$, based on the fact that $V_1$ was maximal, each vertex is either not adjacent to $v_1$ or adjacent to it with both a positive edge and a negative edge. Taking a maximal connected balanced set which does not consists of a single vertex, and applying the same process on $\hat{G_{1}}$ we build a signed graph $\hat{G_{2}}$ with no negative loop where all connections to two of the vertices, when there is any, is by a digon. We may repeat this process until we get a signed graph $\hat{G}^{*}$ where all maximal connected balanced sets are formed of singletons. We may then take $\tilde{H}=\hat{G}^{*}$.
	\end{proof}

	\section{Signed Hadwiger Conjecture}\label{sec:ourHad}
	
	Based on the notion of balanced coloring defined above, we now propose a conjecture that, as we will prove, captures Hadwiger's conjecture and is strongly related to the odd Hadwiger Conjecture, playing an intermediary role between these well-known conjectures.
	
	\begin{conjecture}[Signed-Hadwiger]
		Every signed graph $\hat{G}$ with $\chi_b(\hat{G})\geq t$ has a $\tilde{K}_t$-minor. 
	\end{conjecture}
	
	In this section we discuss the relations between the three different versions of Hadwiger's conjecture. For this we define $\mathcal{SC}_{t}$ to be the class of signed graphs $\hat{G}$ with no negative loop where $\chi_b(\hat{G})\geq t$. Let $SM(t)$ be the largest $k$ such that $\tilde{K}_k$ is a minor of every element of $\mathcal{SC}_{t}$. In this language, the Signed Hadwiger Conjecture claims that $SM(t)=t$.
	
	\subsection{Relating different versions of Hadwiger's conjecture}

	The following theorem can be regarded as a strengthening of a recent result of  Steiner~\cite{S22}, obtained through the notion of balanced coloring.
	
	\begin{theorem}\label{thm:ComparingThreeMinorParameters}
		For every $t$, $M(t)=SM(t)$ and $OM(t)\leq M(t)\leq OM(2t)$.
	\end{theorem}
	
	\begin{proof}
		That $M(t) \geq SM(t)$ follows from the fact that $\tilde{K}_k$ is a minor of $\tilde{G}$ if and only if $K_k$ is a minor of $G$. Similarly, $OM(t)\le M(t)$ follows from the fact that if $G$ has no $K_k$-minor, then $(G, -)$ has no $(K_k,-)$-minor.
		
		To see that $M(t) \leq SM(t)$, let $\hat{G}$ be a signed graph with no $\tilde{K_k}$-minor. Let $\tilde{H}$ be a signed graph obtained from applying Theorem~\ref{thm:MinorAndHomImage}. Since $\tilde{H}$ is a homomorphic image of $\hat{G}$, we have $\chi_b (\tilde{H})\leq \chi_b(\hat{G})$. But $\chi_b (\tilde{H})=\chi(H)$ where $H$ is the underlying graph of $\tilde{H}$ with no loop. On the other hand $\tilde{H}$ is a minor of $\hat{G}$, thus if $\tilde{H}$ contains a $\tilde{K_k}$-minor, then so does $\hat{G}$. However, since each edge of $\tilde{H}$ which is not a loop is part of a digon, if $K_k$ is minor of $H$, then $\tilde{K}_k$ is a minor of $\tilde{H}$, proving that that $M(t) \leq SM(t)$.
		
		For $M(t)\leq OM(2t)$, assuming that every $t$-chromatic graph has a $K_{M(t)}$-minor, we need to show that for every $2t$-chromatic graph $G$, the signed graph $(G,-)$ has a $(K_{M(t)},-)$-minor. First we note that, by the first part of the theorem, our assumption is equivalent to saying that every signed graph of balanced chromatic number $t$ contains a $\tilde{K}_{_{M(t)}}$-minor. Since the only balanced sets in $(G, -)$ are the induced bipartite subgraphs, that $\chi(G)=2t$ implies $\chi_b(G,-)=t$. Then, by the first part of the theorem, $(G,-)$ must contain $\tilde{K}_{_{M(t)}}$-minor which, obviously, contains $(K_{M(t)},-)$ as a subgraph.  
	\end{proof}
	
	We may rephrase this theorem as follows. Assume every $K_t$-minor-free graph (with no loop) is $f(t)$-colorable. Observe that Hadwiger's conjecture claims $f(t)=t-1$, and the best result currently known is that of Delcourt and Postle~\cite{DP21} who show that $f(t)\leq O(t\log \log t)$. The first part of the theorem claims that then every signed graph with no negative loop and no $\tilde{K}_t$-minor admits a balanced $f(t)$-coloring. By associating signed graph $\tilde{G}$ to a graph $G$, we have observed that this statement already contained the original. 
	
	The second part of the theorem is basically the following immediate observation. Again assuming that every $K_t$-minor-free graph (with no loop) is $f(t)$-colorable: if $(G, \sigma)$ has no $(K_t,-)$-minor (hence no $\tilde{K}_t$-minor), then $(G, \sigma)$ admits a balanced $f(t)$-coloring. Observe that in the case of signed graphs of the form $(G,-)$, balanced $f(t)$-coloring is the same as proper $2f(t)$-coloring of $G$; this special case is the result of Steiner~\cite{S22}. The restriction of Theorem~\ref{thm:ComparingThreeMinorParameters} to graphs (no sign) is stated in Theorem~\ref{thm:oddevenHadwiegr}.

	\subsection{Restriction to (unsigned) graphs}\label{sec:unsigned}

	Here, following most of the literature on the odd Hadwiger Conjecture, we avoid signed graphs and define an odd-$K_t$ minor in a graph $G$ as: a 2-coloring of vertices together with a collection $T_1, T_2, ..., T_t$ of vertex disjoint trees in $G$ such that: (i) each edge of any $T_i$ is incident with both colors, and; (ii) between any pair $T_i, T_j$ of trees ($i\neq j$) there is an edge connecting them whose endpoint are of the same color.
	
	Let $f(t)$ be the smallest integer such that every $K_t$-minor-free graph is $f(t)$-colorable. Let $f_o(t)$ be the smallest integer such that every odd-$K_t$-minor-free graph is $f_o(t)$-colorable. Then we can restate Hadwiger's conjecture and the odd Hadwiger Conjecture equivalently as follows. 
	
	\begin{conjecture}[Hadwiger's conjecture, restated] $f(t)=t-1$.
	\end{conjecture}
	
	\begin{conjecture}[Odd Hadwiger Conjecture, restated] $f_o(t)=t-1$.
	\end{conjecture}
	
	The afore-mentioned theorem of Steiner can be stated as follows.
	
	\begin{theorem}[Steiner, \cite{S22}]\label{steiner} $f_o(t)\leq 2f(t)$.
	\end{theorem}
	
	We now introduce the following notion to strengthen Theorem \ref{steiner}.
	
	\begin{definition}
		An \emph{even-odd-$K_t$-minor} of a given graph $G$ is a 2-coloring of vertices together with a collection $T_1, T_2, ..., T_t$ of vertex disjoint trees in $G$ such that: $(i)$ each edge of any $T_i$ is properly colored, and; $(ii)$ between any pair $T_i, T_j$ of trees ($i\neq j$) there is at least one monochromatic edge and at least one properly colored edge.
	\end{definition}
	
	Let $f_{eo}(t)$ be the smallest integer such that every even-odd-$K_t$-minor-free graph is $f_{eo}(t)$-colorable. One may observe that $K_{2t-2}$ has no even-odd-$K_t$-minor, because otherwise two of the trees each have at most 1 vertex and hence the second condition cannot be satisfied between these two trees. Thus $ f_{eo}(t)\geq 2t-2$. However, from Theorem~\ref{thm:ComparingThreeMinorParameters} we get the following.
	
	\begin{theorem}\label{thm:oddevenHadwiegr}
		For every $t\geq 2$ we have $f_o(t)\leq f_{eo}(t)\leq 2f(t)$.
	\end{theorem}

	\section{Topological minors in signed graphs}\label{sec:sudsigned}

	In order to consider subdivisions in signed graphs, we now introduce two definitions which extend the notions of odd-$K_4$ and totally odd-$K_4$, respectively.
	
	\begin{definition}
		A signed graph $(H, \pi)$ is said to be a topological minor of a signed graph $(G, \sigma)$ if: (i) a subdivision of $H$ is isomorphic to a subgraph $G_1$ of $G$, and; (ii)  given any cycle $C$ of $(H, \pi)$ the image of it in $G_1$ has the same sign in $(G, \sigma)$ as the sign of $C$ in $(H,\pi)$.
	\end{definition}
	
	\begin{definition}
		A signed graph $(H, \pi)$ is said to be a total topological minor of a signed graph $(G, \sigma)$ if: (i) a subdivision of $H$ is isomorphic to a subgraph $G_1$ of $G$, and; (ii) given any edge $e$ of $(H, \pi)$, the path $P_e$ representing $e$ in $G_1$ is of the sign $\pi(e)$. 
	\end{definition}

	It follows from the definition that the notion of topological minor is independent of switching. In contrast, the notion of total topological minor is usually based on the choice of the signature. However, there are exceptions and in particular we have the following.
	
	\begin{observation}
		Given a graph $H$, the signed graph $\tilde{H}$ is a topological minor of a signed graph $(G,\sigma)$ if and only if it is a total topological minor of $(G, \sigma)$. 
	\end{observation}
	
	\begin{proof}
		Note that given adjacent vertices $x$ and $y$, in $H$, we have both a negative edge $e^{-}=xy$ and a positive edge $e^{+}=xy$ in $\tilde{H}$. As $\{e^{-},e^{+}\}$ induces a negative 2-cycle in $\tilde{H}$, the paths $P_ {e^{-}}$ and $P_ {e^{+}}$ should be of different signs in $(G, \sigma)$. For each connected pair $xy$ we then associate $e^{-}$ with the negative one of these two and $e^{+}$ with the positive one.
	\end{proof}
	
	Recall that $\mathcal{C}_t$, $\mathcal{OC}_t$, and $\mathcal{SC}_t$ are, respectively: the class of graphs having chromatic number at least $t$; signed graphs of the form $(G,-)$ with $\chi(G)\geq t$, and; signed graphs having balanced chromatic number at least $t$. Based on these notions, we have the following variations of $T(t)$.
	
	\begin{definition} Given a positive integer $t$ we define $OT(t)$ to be the largest $k$ such that $(K_k, -)$ is a topological minor of every member of $\mathcal{OC}_t$. Similarly, $TT(t)$ is the largest $k$ such that $\tilde{K}_k$ is a topological minor of every member of $\mathcal{SC}_t$.
	\end{definition}
	
	\begin{observation}\label{topgens}
		We have $TT(t)\leq OT(t) \leq T(t)$.
	\end{observation}
	
	\begin{proof}
		If $K_k$ is not a topological minor of $G$, then $(K_k,-)$ is certainly not a topological minor of $(G, -)$. And, similarly, if $(K_k, -)$ is not a topological minor of $(G, -)$ then neither is $\tilde{K}_k$. 
	\end{proof}

    In most literature the function $OT(t)$ is studied in its inverse form. That is to ask, given a positive integer $k$, for the largest value of $t$ such that every graph of chromatic number $k$ has a totally odd $K_t$ subdivision, i.e., a subdivision of $K_t$ where each edge is replaced with a path of odd length. The best bound result so far is that of Kawarabayashi~\cite{K2013} which, in our formulation, states:
    
    \begin{theorem}\label{thm:TotallyOddMinor}
    	For any positive integer $t$ we have $OT(t)\geq \sqrt{\frac{4t}{79}}$. 
    \end{theorem}

	{In the next section then we prove the following.}
	
	\begin{theorem}\label{thm:TopMinor-BalancedColoring}
		For any positive integer $t$ we have $TT(t)\geq \sqrt{\frac{2t}{79}}$. 
	\end{theorem}
	
	Restricted to (unsigned) graphs this can be restated as: any graph of chromatic number high enough, not only contains a totally odd $K_t$ subdivision, but also contains a set of~$t$ vertices where between each pair of them there is a path of odd length and path of even length where no pair of these paths have a common internal vertex.

	\subsection{Topological minors and balanced coloring}
	
	Here we show the connection between the absence of a large topological minor and balanced coloring in signed graphs by proving the following stronger version of  Theorem~\ref{thm:TopMinor-BalancedColoring}.
	
	\begin{theorem}\label{thm:Signed Subdivision}
		Let $G=(V,E)$ be a signed graph with no $\tilde{K}_t$-subdivision. For any vertex set $Z\subseteq V$ with $|Z| \leq 2t^2$ any precoloring of the subgraph of $G$ induced by $Z$ can be extended to a $\frac{79}{2}t^2$-coloring of $G$.
	\end{theorem}
	
	The proof given below is an adaptation of the proof by Kawarabayashi for the existence of large totally odd subdivisions in graphs of high chromatic number \cite{K2013}. We first state some key results from the literature that are needed for the proof. The first one is the following folklore observation.
	
	\begin{observation}\label{obs:spanning balanced subgraph}
		Let $\hat{G}$ be a signed graph and assume $\hat{H}$ is a balanced subgraph of $\hat{G}$ with the maximum possible number of edges. Then for each vertex $v$ of $G$ we have $d_{H}(v)\geq \frac{d_{G}(v)}{2}$. In particular $\delta(H)\geq \frac{\delta(G)}{2}$. 
	\end{observation}
	
	For the proof the observation one may note that $H$ can be regarded as the spanning subgraph  of the positive edges of $(G,\sigma')$ for some switching $\sigma'$ of $\sigma$. Now, if a vertex $v$ is incident to less than $\frac{d_{G}(v)}{2}$ positive edges in $(G,\sigma')$, then after a switching at $v$ we have a signature $\sigma''$ with more positive edges, contradicting the choices of $H$. 
	
	The next statement is obtained from \cite{CGGGLS06} by taking signed graphs as symmetric group labeled digraphs, with the group being (additive) $\mathbb Z_{2}$ where $0$ plays the role of $+$ and $1$ plays the role of $-$.  
	
	\begin{theorem}\label{thm:negative paths}
		Let $G$ be a signed graph and $H$ be a balanced, connected subgraph so that all edges of $H$ are positive. For any fixed $k$ one of the following holds.
		\begin{enumerate}
			\item There are $k$ mutually disjoint negative $H$-paths, i.e., $k$ mutually disjoint paths $P_1, \dots, P_k$ such that each $P_i$ is a negative path whose end vertices are in $H$, or
			\item There is a set $X \subseteq V(G)$ of at most  $2k-2$ vertices, such that every negative path with both end vertices in $H$ contains a vertex in $X$.
		\end{enumerate}
	\end{theorem}

	A graph is \emph{$\ell$-linked} if it has at least $2\ell$ vertices, and for any choice of distinct vertices $u_1,u_2, \dots, u_\ell,\allowbreak v_1, v_2, \dots, v_\ell$ there are $\ell$ mutually disjoint paths $P_1, P_2, \dots P_\ell$ so that $P_i$ has ends $u_i,v_i$ ($1 \leq i \leq \ell$). The following statement is from \cite{TW2005}.
	
	\begin{theorem}\label{thm:l-linked}
		Every $2\ell$-connected graph $G$ with at least $5\ell|V(G)|$ edges is $\ell$-linked.
	\end{theorem}
	This theorem will be used in combination with the following result from \cite{BKMM2009}.
	\begin{theorem}\label{thm:l-linkedcomb}
		Let $G$ be a graph and $k$ an integer such that
		$$|V (G)|\geq \frac{5}{2} k \,\,\,\,\, \text{and} \,\,\,\,\,
		|E(G)| \geq \frac{25}{4} k |V (G)| - \frac{25}{2} k^2.$$
		Then $|V (G)| \geq 10k + 2$ and $G$ contains a $2k$-connected subgraph $H$ with at least $5k|V (H)|$ edges.
	\end{theorem}

	\begin{proof}[Proof of Theorem~\ref{thm:Signed Subdivision}]
		
		Let $\hat{G}$ be a minimum counterexample with respect to the number of vertices. That is, $\hat{G}$ does not have a $\tilde{K}_t$-subdivision and there exists $Z \subseteq V$ with $|Z| \leq 2t^2$ and a precoloring of $Z$ that cannot be extended to a $\frac{79}{2}t^2$-coloring of $G$. Then $\hat{G}$ must have at least $\frac{79}{2}t^2 +1$ vertices. We prove several claims about $\hat{G}$ before getting to a contradiction by providing a $\tilde{K}_t$-subdivision.
		
		\newcounter{count}
		\refstepcounter{count}\label{count: no parallel}
		(\thecount) \emph{We may assume $\hat{G}$ has no parallel edges of the same sign.}
		
		This is because such parallel edges do not affect any coloring, and if $\hat{G}$ already has no $\tilde{K}_t$-subdivision, then it certainly does not have one after deleting an edge. 
		
		\refstepcounter{count}\label{count: min degree}
		(\thecount) \emph{Every vertex $v \in V-Z$ has degree at least $\frac{79}{2}t^2$ in $\hat{G}$.}
		
		Otherwise, by minimality, for some $v$ with degree at most $\frac{79}{2}t^2-1$, the signed graph $\hat{G}-v$ has a $\frac{79}{2}t^2$-coloring which is an extension of the precoloring of $Z$. But because $v$ has low degree, this coloring of $\hat{G}-v$ can be extended to $\hat{G}$, a contradiction.

		An \emph{$\ell$-separation} of $\hat{G}$ is a pair $(\hat{G}_1,\hat{G}_2)$ of subgraphs so that $\hat{G}_1 \cup \hat{G}_2 =\hat{G}$, and $|V(\hat{G}_1) \cap V(\hat{G}_2)| = \ell$. Following Kawarabayashi, we say that an $\ell$-separation $(\hat{G}_1,\hat{G}_2)$ is \emph{$Z$-essential} if each $\hat{G}_i$ ($i=1,2$) has at least one vertex which is not in $\hat{G}_j\cup Z$ for $j\neq i$.  
		
		\refstepcounter{count}\label{count:no Zessential sep}
		(\thecount) \emph{For $\ell \leq t^2$, $\hat{G}$ admits no $Z$-essential $\ell$-separation.}
		
		Suppose for a contradiction that such a separation $(\hat{G}_1,\hat{G}_2)$ exists. Consider the $V(\hat{G}_1)\setminus {V(\hat{G}_2)}$ or $V(\hat{G}_2) \setminus V(\hat{G}_1)$ which has no common element. Since $|Z| \le 2t^2$ the number of elements of $Z$ in one of these two sets is at most $t^2$. By symmetry, we may assume that $|Z \cap V(\hat{G}_1)\setminus {V(\hat{G}_2)} |\leq t^2$.
		
		By the minimality of $\hat{G}$ and since there is at least one vertex of $\hat{G}$ not in $V(\hat{G}_2)\cup Z$, the precoloring $\varphi$ of $Z$ can be extended to a coloring $\varphi'$ of $\hat{G}_2\cup Z$. Now consider the restriction $\varphi'$ of $\varphi$ on the vertices of $\hat{G}_1$ that are colored. Observe that there are at most $2t^2$ such vertices and that $\hat{G}_1$ has at least one less vertex than $\hat{G}$. Thus, by the assumption on the minimality of $\hat{G}$, the coloring $\varphi'$ can be extended to the rest of $\hat{G}_1$, resulting a coloring of $\hat{G}$, a contradiction.

		\refstepcounter{count}\label{count: span bal subgraph}
		(\thecount) \emph{There is a spanning balanced subgraph $\hat{H}$ of $\hat{G}-Z$ whose minimum degree is at least $\frac{75}{4}t^2$.}

		It follows from (\ref{count: min degree}) that $\delta(\hat{G}-Z) \geq \frac{79}{2}t^2-2t^2=\frac{75}{2}t^2$. The claim then follows by Observation~\ref{obs:spanning balanced subgraph}. 
		
		In the rest of the proof we will assume the signature of $\hat{G}$ is switched, if needed, so that all edges of $H$ are positive.

		\refstepcounter{count}\label{count: linked subgraph}
		(\thecount) \emph{There is a subgraph $L \subseteq H$ which is $\frac{3}{2}t^2$-linked, $3t^2$-connected, and, in particular, has at least $3t^2$ vertices.}
		
		From (\ref{count: no parallel}), and because it is balanced, $H$ has no parallel edges or digons. From (\ref{count: span bal subgraph}), $H$ has minimum degree at least $\frac{75}{4}t^2$, and because these neighbors are distinct, $H$ has at least these many vertices and at least $\frac{75t^2}{8}|V(H)|$ edges.
		We may then apply Theorem \ref{thm:l-linkedcomb} with $k=\frac{3t^2}{2}$ to get a subgraph $L$ of $H$ which is $3t^2$-connected and with at least $15t^2 |V(L)|$ edges. Now, taking $\ell=\frac{3}{2}t^2$, Theorem \ref{thm:l-linked} ensures that $L$ is $\frac{3}{2}t^2$-linked.

		Recall that, being a subgraph of $H$, all edges of $L$ are positive in the signature of $\hat{G}$ that we are working with.
		
		\refstepcounter{count}\label{count: par break paths}
		(\thecount) \emph{There are $\frac{1}{2}t^2 \geq \binom{t}{2}$ disjoint negative $L$-paths in $G$.}
		
		Suppose for a contradiction that such paths do not exist. Then by Theorem~\ref{thm:negative paths}, there is a subset $X \subseteq V(G)$ with $|X|\leq t^2-2$ so that $G-X$ has no negative $L$-path. From (\ref{count: linked subgraph}), $L-X$ is $2$-connected, and it is, therefore, contained in some $2$-connected block $L'$ of $G-X$. We now prove two claims in order to prove \eqref{count: par break paths}.
		
		(\thecount a) \emph{$L'$ is balanced.} If not, then there is a negative cycle $C \subseteq L'$. Then due to the $2$-connectivity of $L'$, there exist two disjoint paths (possibly trivial ones) in $L'$, joining $C$ and $L$. However, this structure contains a negative $L$-path, a contradiction. 
		
		If $L'=G-X$, then we can extend a precoloring of $Z$ to a  $3t^2$-coloring of $G$ as follows: the precoloring of $Z$ uses at most ${2t^2}$ colors, then at most a set of $t^2-2$ colors are used for coloring vertices in $X-Z$. Finally one new color is needed for all vertices in $L'$ due to  (\ref{count: par break paths}a). This is a contradiction as $3t^2 < \frac{79}{2}t^2$.

		Let $W_1, W_2, \dots W_r$ be the remaining $2$-connected blocks in $G-X$ for $r \geq 1$. Denote by $v_i$, the cut-vertex in $V(L') \cap V(W_i)$, if one exists.

		(\thecount b) \emph{$W_i - v_i \subseteq Z$ for $1 \le i \le r$.} Observe that, $|L'| \ge |L-X| \ge |Z| + 2$, where the second inequality follows from (\ref{count: linked subgraph}) because $|X| \le t^2-2$ and $|Z| \leq 2t^2$. So if there is a $v \in W_i - v_i$ such that $v \notin Z$, then there is a $Z$-essential separation of order at most $|X| + 1 \le t^2 - 1$, contradicting (\ref{count:no Zessential sep}).  
		
		Now we can extend the precoloring of $Z$ to a $\frac{79}{2}t^2$-coloring of $G$ as before: use at most $2t^2$ colors in the precoloring of $Z$, observe that $G-Z \subseteq L' \cup X$ by (\ref{count: par break paths}b), then use at most $t^2-2$ additional colors to color the remainder of $X$ because $|X| \le t^2-2$, and
		use one additional color to color the remainder of $L'$ by (\ref{count: par break paths}a).
		
		This coloring uses at most $3t^2$ colors, which contradicts that $G$ is a counterexample to the theorem. This completes the proof of (\ref{count: par break paths}).
		
		Now we will demonstrate that there is a $\tilde{K}_t$-subdivision in $G$. We will construct this subdivision using $L$ and negative $L$-paths. From (\ref{count: par break paths}), there exist  $\frac{1}{2}t^2 \geq \binom{t}{2}$ disjoint negative $L$-paths in $G$. Choose $\binom{t}{2}$ such paths and let $W$ denote their endpoints. We have $|W|=t(t-1) < t^2$. By (\ref{count: linked subgraph}), $L$ has at least $3t^2$ vertices, and so we may choose $t$ distinct vertices $u_1, \dots, u_t$ in $L-W$. These will serve as the terminals of the $\tilde{K}_t$-subdivision. 
		
		For each pair of distinct terminals $u_i,u_j$, (there are exactly $\binom{t}{2}$ such pairs), 
		we associate exactly one negative $L$-path, $P_{ij}$. Refer to the ends of $P_{ij}$ as $p_{ij}$ and $p_{ij}'$. Furthermore, for each vertex $u_i$ choose a set of neighbours $N_i$ in $L-W$ of size $2(t-1)$, say  
		$$N_i= \{ w_{i,1}, \ldots, w_{i,i-1}, w_{i,i+1},\ldots w_{i,t}, v_{i,1}, \ldots,  v_{i,i-1}, v_{i,i+1},\ldots v_{i,t}\}$$
		
		such that $\{N_i\cup u_i\} \cap  \{N_j\cup u_j\}=  \emptyset$ for each $j \neq i$. It is possible to do so since by (\ref{count: linked subgraph}), $L$ is $3t^2$-connected, and hence, the minimum degree of $L$ is at least  $3t^2$ which is bigger than 
		$$ 2t(t-1)+2 \binom{t}{2}+ t.$$
		
		Now, we find the following disjoint paths in $L$.
		\begin{enumerate}
			\item For each pair $i, j $ in $\binom{[t]}{2}$, a path with ends $w_{i,j} w_{j,i}$. This will serve as the positive paths in the $\tilde{K}_t$-subdivision.
			\item For each pair $i, j $ in $\binom{[t]}{2}$, two paths, one with ends $v_{i,j},p_{ij}$ and one with ends $v_{j,i}, p_{ij}'$, for $i<j$. These paths together with $P_{ij}$ will serve as the negative paths in the $\tilde{K}_t$-subdivision.
		\end{enumerate}

		This is a total of $3 \binom{t}{2}$ disjoint paths in $L$. Since, by (\ref{count: linked subgraph}), $L$ is $\frac{3}{2}t^2$-linked, and $3 \binom{t}{2} \leq \frac{3}{2}t^2$,
		we will be able to do so.

		This means there is a $\tilde{K}_t$-subdivision in $G$, which contradicts our choice of $G$ and completes the proof.
	\end{proof}

	\section{Fractional signed Hadwiger}
	
	The fractional balanced chromatic number of a signed graph with no negative loop can be defined analogous to the classic fractional chromatic number: A \emph{fractional balanced coloring} of a signed graph $(G,\sigma)$ is to assign non-negative weights to the balanced sets of a $(G,\sigma)$ such that each vertex is in balanced sets of total weight at least 1. The \emph{fractional balanced chromatic number} of $(G, \sigma)$, denoted $\chi_{fb}(G,\sigma)$, is the infimum of total weights among all fractional balanced colorings. For a comprehensive study of fractional balanced coloring we refer to \cite{KNWYZZ25+}. Here as a further application of Theorem~\ref{thm:MinorAndHomImage} we state the following.
	
	\begin{theorem}
		Given a signed graph $(G, \sigma)$ with no $\tilde{K}_t$-minor, we have $\chi_{fb}(G,\sigma)\leq 2t-2$. 
	\end{theorem}	
	
	This is a corollary of the following theorem from \cite{RS98} and the observation that $\chi_{fb}(\tilde{H})=\chi_f(H)$, noting that positive loops have no affect on balanced fractional coloring of a signed graph. We leave the detail to the reader.
	
	\begin{theorem}\cite{RS98}
		Given a graph $G$ with no $K_t$-minor we have $\chi_f(G)\leq 2t-2$.
	\end{theorem}

	\section{Concluding remarks}
	
	In this work we introduced a signed version of Hadwiger's conjecture and showed that while it helps to bound the chromatic number of dense families of (signed) graphs it is still equivalent to the Hadwiger's conjecture which only bounds the chromatic number of sparse families of graphs. Our conjecture in turn helps to better understand the connection between Hadwiger's conjecture and the odd Hadwiger Conjecture.
	
	A natural line of work to improve on the existing bounds would be to consider signed simple graphs, or more generally signed graphs of given girth. For signed planar simple graphs the best upper bound for the circular chromatic number is 6, while a construction for a simple planar graph of circular chromatic number $\frac{14}{3}$ is given in \cite{NWZ21}. The exact value remains open. There are no specific improvements on the corresponding bounds for other classes of signed graphs such as signed $K_{t}$-minor-free simple graphs.
	
		We also showed the existence of relatively large subdivisions in signed graphs of high balanced chromatic number. When restricted on graphs this shows the existence of a subdivision where between any pair of vertices there are disjoint odd and even paths assuming the chromatic number is high. While we expect that the bounds given here can be improved, we also believe the bounds for chromatic number to enforce such stronger structure would be higher. Indeed to have both odd and even path between $t$ vertices one needs at least ${t+1 \choose 2}$ vertices which is already much higher lower bound. It is also not known how assumptions such as high girth would affect this bound. 


    { While this work was under review, M. Kühn, L.  Sauermann, R. Steiner and Y. Wigderson \cite{KSSW2025+} announced a  disprove of the Odd Hadwiger conjecture, showing that for $n$ large enough there are graphs on $n$ vertices of independence number 2 (thus chromatic number at least $\frac{n}{2}$) whose largest odd $K_t$-minor is of order not much exceeding $\frac{n}{3}$. This does not disprove Hadwiger's conjecture. }

	{\bf Acknowledgment.} This work is supported by the following grants and projects: 1. ANR-France project HOSIGRA (ANR-17-CE40-0022). 2. MATH-AMSUD MATH230035 and MATH210008. 3. ANID/Fondecyt Regular 1220071 and ANID-MILENIO-NCN2024-103 (Andrea Jim\'enez). 4. Simons Foundation grant \#845698 (Jessica McDonald) 5. NSERC (Kathryn Nurse) 6. ANID/Fondecyt Regular 1252197 (Daniel A. Quiroz).

\end{document}